\documentclass[reqno]{pspum-l}
\usepackage{amsmath,amsthm,amscd,amssymb}
\usepackage{latexsym}

\usepackage{hyperref}

\issueinfo{87}
  {0}
  {Xxxx}
  {2013}
\pagespan{343}{364}
\copyrightinfo{2012}{(The authors)}

\newcommand*{\mailto}[1]{\href{mailto:#1}{\nolinkurl{#1}}}
\newcommand{\arxiv}[1]{\href{http://arxiv.org/abs/#1}{arXiv:#1}}

\newtheorem{theorem}{Theorem}[section]

\theoremstyle{definition}

\newcommand{\R}{{\mathbb R}}
\newcommand{\N}{{\mathbb N}}

\newcommand{\C}{{\mathbb C}}

\newcommand{\be}{\begin{equation}}
\newcommand{\ee}{\end{equation}}

\newcommand{\ol}{\overline}
\newcommand{\ti}{\tilde}

\newcommand{\bay}{\begin{array}}
\newcommand{\eay}{\end{array}}

\newcommand{\spr}[2]{\langle #1 , #2 \rangle}

\newcommand{\E}{\mathrm{e}}
\newcommand{\I}{\mathrm{i}}

\newcommand{\im}{\mathrm{Im}}
\newcommand{\re}{\mathrm{Re}}
\newcommand{\db}{\mathrm{dom}}
\newcommand{\hr}{\mathcal{H}}
\DeclareMathOperator{\ran}{ran}

\newcommand{\eps}{\varepsilon}

\newcommand{\sig}{\sigma}
\newcommand{\lam}{\lambda}
\newcommand{\gam}{\gamma}

\numberwithin{equation}{section}

\begin{document}

\title{Spectral Theory as Influenced by Fritz Gesztesy}

\author[G.\ Teschl]{Gerald Teschl}
\address{Faculty of Mathematics\\ University of Vienna\\
Nordbergstrasse 15\\ 1090 Wien\\ Austria\\ and International
Erwin Schr\"odinger
Institute for Mathematical Physics\\ Boltzmanngasse 9\\ 1090 Wien\\ Austria}
\email{\mailto{Gerald.Teschl@univie.ac.at}}
\urladdr{\url{http://www.mat.univie.ac.at/~gerald/}}

\author[K.\ Unterkofler]{Karl Unterkofler}
\address{Vorarlberg University of Applied Sciences\\
Hochschulstr.~1\\
6850 Dornbirn\\ Austria}
\email{\mailto{Karl.Unterkofler@fhv.at}}
\urladdr{\url{https://homepages.fhv.at/ku/}}

\dedicatory{To Fritz Gesztesy, teacher, mentor, and friend, on the occasion of his 60th birthday.}
\thanks{in "Spectral Analysis, Differential Equations and Mathematical Physics", 
H. Holden et al. (eds), 343--364, Proceedings of Symposia in Pure Mathematics {\bf 87}, Amer. Math. Soc., Providence, 2013.}
\thanks{Research supported by the Austrian Science Fund (FWF)  under Grant No.\ Y330 and
P24736-B23}


\begin{abstract}
We survey a selection of  Fritz's principal contributions to the field of spectral theory and, in particular, to Schr\"odinger operators.
\end{abstract}

\maketitle

\section{Introduction}

The purpose of this Festschrift contribution is to highlight some of Fritz's profound contributions to spectral theory and, in
particular, to Schr\"odinger operators. Of course, if you look at his list of publications it is clear that this is a {\em mission impossible}
and hence the present review will only focus on a small selection. Moreover, this selection is highly subjective and biased by
our personal research interests:

\begin{itemize}
\item Relativistic Corrections
\item Singular Weyl--Titchmarsh--Kodaira Theory
\item Inverse Spectral Theory and Trace Formulas
\item Commutation Methods
\item Oscillation Theory
\item Non-self-adjoint operators
\end{itemize}

Other aspects of his work, e.g., on point interactions and integrable nonlinear wave equations are summarized in the monographs \cite{aghh} and \cite{ghV1,ghmtV2},
respectively.
For some of his seminal contributions to inverse spectral theory \cite{drgs,gsipa1,gstsp,gsipa2,gslbm,gsun,gsna} we refer to Fritz's own review \cite{gefest}.
For some of his recent work on partial differential operators, Krein-type resolvent formulas, operator-valued Weyl--Titchmarsh operators (i.e., energy-dependent analogs of Dirichlet-to-Neumann maps), and Weyl-type spectral asymptotics for Krein--von Neumann extensions of the Laplacian on bounded domains, we refer to \cite{AGMT10,GKMT01,GM08,GM11,GMZ07}, and the detailed lists of references therein.

To be sure, this is not a survey of the state of the art on these topics but we rather focus exclusively on Fritz Gesztesy's contributions to and influence exerted on the field. Especially,
the bibliography is far from being complete and only reflects the particular purpose of this survey.

{\bf Personal note.} It was 28 years ago in spring of 1985 that Fritz and Karl first met in Graz. 
Fritz, jointly with Wolfgang Bulla, advised his master and PhD thesis.
Karl then spent two years (1990--1992) as a postdoc with Fritz at the University of  Missouri, Columbia.
When he returned to Austria, Karl met Gerald, who  had just finished his masters with Wolfgang Bulla.
Karl immediately proposed to him to pursue his PhD with Fritz. After a brief meeting with Fritz in Graz
and a short visit to the University of  Missouri, Columbia, Gerald enrolled in the PhD program there
and a second successful cooperation began. Fritz guided us through the early stages of our careers and was much more than a mentor for us. On the one hand he
was an unlimited source of ideas for new projects and on the other hand he was available for discussing mathematics
close to 24 hours at 358 days a year (one week has to be subtracted for which he disappears when snorkeling in the Caribbean or Hawaii).
His devotion to mathematics, his perfection when it comes to details, and his endurance has always been a role model for us.
And last but not least, he is the one everybody asks for references; his collection of papers and books is legendary!

\begin{center}
Happy Birthday, Fritz, and on to many more such anniversaries!
\end{center}

\section{Relativistic Corrections}

In 1950 Foldy and Wouthuysen developed a formal  perturbation scheme 
(which is now ``regarded as obsolete'', according to Kutzelnigg on page 678 in \cite{Schwerdtfeger}) to obtain relativistic corrections of the nonrelativistic Pauli operator $H_{+}$. However,  adding the first
correction term of order $c^{-2}$  already destroys all spectral properties of $H_{+}$. 
The fact, that nevertheless formal perturbation theory  yields correct results in   special
cases has been explained in terms of spectral concentration  in \cite{Gesztesy85,Gesztesy82} and \cite{Veselic69}.

Historically, the first rigorous treatment of the nonrelativistic limit 
of Dirac Hamiltonians  goes back to Titchmarsh \cite{Titchmarsh62} who 
proved holomorphy of the Dirac eigenvalues (rest energy subtracted) with
respect to $c^{-2}$ for spherically symmetric potentials and obtained explicit
 formulas for relativistic bound state corrections of order O$(c^{-2})$
(formally derived in \cite{Sewell49}). Holomorphy of the Dirac resolvent in three
dimensions in $c^{-1}$ for electrostatic interactions were first obtained
by Veselic \cite{Veselic71} and then extended to electromagnetic interactions by 
Hunziker \cite{Hunziker75}. An entirely different approach, based on an abstract set up,
has been used in \cite{Cirincione81} to prove strong convergence of the unitary 
groups as $c^{-1} \to 0$.  Convergence of solutions of the Dirac equation 
based on semigroup methods has also been obtained in \cite{Schoene79}.
Time-dependent electromagnetic fields are treated in \cite{mauser}.

In joint collaboration with Bernd Thaller and Harald Grosse \cite{Gesztesy84} (see also \cite{Gesztesy83} and \cite{Gesztesy84b})
employing an abstract framework, holomorphy of the Dirac resolvent in $c^{-1}$ under very general conditions on the electromagnetic interaction potentials was obtained. 
Moreover, this approach led to the first rigorous
derivation of explicit formulas for relativistic corrections of order
O$(c^{-2})$ to bound state energies.

Application of these results lead to
relativistic corrections for energy bands and corresponding corrections 
for impurity bound states for one-dimensional periodic systems
 \cite{bulla88}   and relativistic corrections for the scattering matrix   \cite{bulla92}. 
An explicit treatment of relativistic corrections of the 
scattering amplitude appeared in \cite{Grigore89}.
 
A thorough and comprehensive treatment of  Dirac operators can be found in the book by Thaller \cite{thallerbook} and in \cite{Schwerdtfeger}. \\

Based on the abstract approach introduced in  \cite{Hunziker75}  and \cite{Cirincione81}, 
let $ \hr_{\pm}   $ be separable, complex Hilbert spaces. 
One introduces self-adjoint operators $\alpha, \beta$ in
$\hr = \hr_{+} \oplus \hr_{-}$ of the type
\begin{align}
\alpha = 
 \begin{pmatrix}
0 & A^{*} \\
A & 0 
 \end{pmatrix} , \quad \ 
\beta =
 \begin{pmatrix}
1 & 0 \\ 
0 & -1 
 \end{pmatrix} ,
\end{align}
where $A$ is a densely defined, closed operator from 
$\hr_{+}$ into $\hr_{-}$. Next,  define
the abstract free Dirac operator $H^{0}(c)$ by
\begin{align}
H^{0}(c) = c\, \alpha + mc^{2} \beta, \quad
 \db(H^{0}(c)) =
\db(\alpha), \quad
c \in \R \setminus \{ 0 \}, \ \ \ m > 0 
\end{align}
and the interaction $V$ by
\begin{align}
V =  
 \begin{pmatrix}
V_{+} & 0 \\
0 & V_{-} 
 \end{pmatrix} ,
\end{align}
where $V_{\pm}$ denote self-adjoint operators in $\hr_{\pm}$,
respectively.
 Assuming $V_{+}$ (resp.$V_{-}$) to be bounded
w.r.t.\  $A$ (resp. $A^{*}$), i.e.,
\begin{align}
 \db(A) \subseteq \db(V_{+}), \quad
 \db(A^{*}) \subseteq \db(V_{-}),
\end{align}
the abstract Dirac operator $H(c)$ reads
\begin{align}\label{diracoperator}
H(c) = H^{0}(c) + V, \quad \db(H(c)) = \db(\alpha).
\end{align}
Obviously, $H(c)$ is self-adjoint for $|c|$ large enough.
The corresponding self-adjoint (free) Pauli operators
in $\hr_{\pm}$ are then defined by
\begin{align}
& H_{+}^{0} = (2m)^{-1} A^{*}A, \quad   H_{+} = H_{+}^{0} + V_{+},  \quad \db (H_{+}) = \db(A^{*}A),  \nonumber \\  
&H_{-}^{0} = (2m)^{-1} AA^{*},  \quad
  H_{-} = H_{-}^{0} + V_{-},   \quad 
 \db (H_{-}) = \db(AA^{*}).
\end{align}
Introducing the operator $B(c)$ (see \cite{Hunziker75})      
\begin{align}
B(c) =
 \begin{pmatrix}
1 & 0 \\
0 & c 
 \end{pmatrix} ,
\end{align}
one establishes the holomorphy of the Dirac resolvent in $c^{-1}$.

\begin{theorem} \label{t21}
 Let $H(c)$ be defined as above and fix 
$z \in \ \C \setminus \R$. Then \\
 (a) $(H(c)-mc^{2} -z)^{-1}$ is holomorphic w.r.t.\   $c^{-1}$ around
$c^{-1}=0$
\begin{align}
&  (H(c)-mc^{2} -z)^{-1}  \nonumber \\
& = \left\{ 1+
 \begin{pmatrix}
0 \ \ &
(2mc)^{-1}(H_{+}-z)^{-1}A^{*}(V_{-}-z) \\
(2mc)^{-1}A(H_{+}^{0}-z)^{-1}V_{+} \ \ & 
(2mc^{2})^{-1}z(H_{-}^{0}-z)^{-1}(V_{-}-z)
 \end{pmatrix}
 \right\}^{-1} \times  \nonumber \\
& \times
 \begin{pmatrix}
(H_{+}-z)^{-1} \ \ &
(2mc)^{-1}(H_{+}-z)^{-1}A^{*} \\
(2mc)^{-1}A(H_{+}^{0}-z)^{-1} \ \  &
(2mc^{2})^{-1}z(H_{-}^{0}-z)^{-1}
 \end{pmatrix}
. \label{110}
\end{align}
(b)
$B(c)(H(c)-mc^{2} -z)^{-1}B(c)^{-1}$ is holomorphic  w.r.t.\
 $c^{-2}$ around
$c^{-2}=0$
\begin{align}
& B(c)(H(c)-mc^{2} -z)^{-1}B(c)^{-1}   \nonumber \\
& = \left\{ 1+
 \begin{pmatrix}
0 \ \ &
(2mc^{2})^{-1}(H_{+}-z)^{-1}A^{*}(V_{-}-z) \\
0  \ \ & 
(2mc^{2})^{-1} [(2m)^{-1} A(H_{+}-z)^{-1}A^{*}-1](V_{-}-z)
 \end{pmatrix}
 \right\}^{-1} \times  \nonumber \\
&  \times 
 \begin{pmatrix}
(H_{+}-z)^{-1} \ \ &
(2mc^{2})^{-1}(H_{+}-z)^{-1}A^{*} \\
(2m)^{-1}A(H_{+}-z)^{-1}  \ \ & 
(2mc^{2})^{-1}[(2m)^{-1}A(H_{+}-z)^{-1}A^{*}-1]
 \end{pmatrix}
 .    \label{111}
\end{align}
First order expansions in \eqref{110} and \eqref{111} yield
\begin{align} \label{210}
&(H(c)-mc^{2} -z)^{-1}  = 
 \begin{pmatrix}
(H_{+}-z)^{-1} \ & 0 \\
0 \ & 0 \\
 \end{pmatrix}
 \nonumber \\ & + c^{-1}
 \begin{pmatrix}
0 \ \ &
(2m)^{-1}(H_{+}-z)^{-1}A^{*} \\
(2m)^{-1}A(H_{+}-z)^{-1}  \ \ &  0
 \end{pmatrix}
   + O(c^{-2})   
\end{align}
(clearly illustrating the nonrelativistic limit $|c| \to \infty$) and
\begin{align}
& B(c)(H(c)-mc^{2} -z)^{-1}B(c)^{-1} =: R_{0}(z) +c^{-2} R_{1}(z) + O(c^{-4})  \nonumber \\
& = 
 \begin{pmatrix}
(H_{+}-z)^{-1} \ \ &  0 \\
(2m)^{-1}A(H_{+}-z)^{-1} \ \ & 0
 \end{pmatrix}
 + c^{-2} 
 \begin{pmatrix}
R_{11}(z) \ \ & R_{12}(z) \\
R_{21}(z) \ \ & R_{22}(z)
 \end{pmatrix}
  + O(c^{-4}) , \nonumber 
  \end{align}
 where
\begin{align}
& R_{11} (z) = (2m)^{-2}  (H_{+} - z)^{-1} A^{*} 
(z-V_{-}) A (H_{+} - z)^{-1},  \nonumber \\
& R_{12} (z) = (2m)^{-1}  (H_{+} - z)^{-1} A^{*},   \nonumber \\
& R_{21} (z) = (2m)^{-2} \left( (2m)^{-1} A (H_{+} - z)^{-1} A^{*} -1 \right)
(z-V_{-}) A (H_{+} - z)^{-1},  \nonumber \\
& R_{22} (z)= (2m)^{-1} \left( (2m)^{-1} A (H_{+} - z)^{-1} A^{*} -1 \right). 
 \label{n214} 
\end{align}
\end{theorem}
To formulate and prove Theorem \ref{t21} Fritz et al.\ combined a clever decomposition of the Dirac resolvent with a skillful application of some commutations formulas by Deift \cite{de} (see also Section~\ref{commutation}). 

Analyzing the relationship between the spectrum of $(H_{ +} - z)^{-1}$ and $R_{0}(z)$ now yields the following result on relativistic eigenvalue corrections.

\begin{theorem}\label{t22}
 Let $H(c)$ be defined as in \eqref{diracoperator} and assume $E_{0} \in \sigma_{d}(H_{ +})$ to be a discrete
eigenvalue of $H_{ +}$  of multiplicity $m_0 \in \N$. Then, for $c^{ - 2}$ small enough, $H (c) - m  c^{2} $ has precisely $m_0$ eigenvalues (counting multiplicity) near $E_0$  which are all holomorphic w.r.t.\   $c^{-2}$.
 More precisely, all eigenvalues $E_j(c^{- 2})$ of $H(c) - m\, c^{2}$ near $E_0$ satisfy
\begin{align} \label{eigenvalues}
E_j(c^{-2}) = E_{0} + \sum_{k=1}^{\infty} (c^{-2})^k  E_{j,k},   \quad  j = 1, \ldots ,j_0, \quad j_0 \leq m_0
\end{align}
and if $m_j$ denotes the multiplicity of $E_j(c^{ -2})$ then $\sum_{j=1}^{j_0}m_j =m_0$.\\
In addition, there exist linearly independent vectors
\begin{align}
f_{jl}(c^{-1}) = 
\begin{pmatrix} f_{+jl} (c^{-2}) \\
c^{-1} f_{-jl}(c^{-2})
\end{pmatrix}, \ \ j= 1, \ldots j_{0}, \ \ l=1,\ldots, m_j
\end{align}
s.t.\ $f_{ \pm ji}$ are holomorphic w.r.t.\  $c^{ -2}$ near $c^{ -2}=0$ and
\begin{align}
& H_{+} f_{+jl}(0) = E_0 f_{+jl}(0) , \quad f_{-jl}(0)= (2 m)^{-1} A  f_{+jl}(0) 
\end{align}
 and
\begin{align}
& (H(c)-m c^{2})f_{jl}(c^{-1}) = E_j(c^{-2}) f_{jl}(c^{-1}) , \quad  j= 1, \ldots ,j_0,
\quad  l= 1, \ldots ,m_j .  
\end{align}
The eigenvectors $f_{jl}(c^{-1})$ can be chosen to be orthonormal. Finally, the first-order corrections
$ E_{j,1}$ in \eqref{eigenvalues} are explicitly given as the eigenvalues of the matrix
\begin{align}
(2m)^{-2}\left( A  f_r,(V_{-}- E_0) A f_s\right), \quad r, s = 1, \ldots\ldots , m_0 ,  
\end{align}
where $\{ f_r \}_{r=1}^{m_0}$ is any orthonormal basis of the eigenspace of $H_{+}$ to the eigenvalue $E_0$.
\end{theorem}
Remarks:  
 (a) The main idea  behind Theorem \ref{t22} was to look for eigenvalues
of the resolvent $(H (c) - m c^{2} - z)^{ -1}$  and applying the strong spectral mapping theorem
(\cite{ReedSimon4},  page 109, Lemma 2) instead of looking directly for eigenvalues of the unbounded Hamiltonian
$ H(c) - m c ^{2}$.

(b) Theorem \ref{t22} for $m_0 = 1$ is due to Fritz et al.\ (see \cite{Gesztesy84,Gesztesy83}). In the general case $m_0 > 1$ only holomorphy of $E_j(c^{-1})$ w.r.t.\  $c^{ -1}$ near $c^{ -1} = 0$ and
$ E_j (c^{ -1}) - m c^{2} - E_0 =_{c \to \infty} O(c^{ -2})$ has been proven in \cite{Gesztesy84}. 
The above extension of this result for $m_0 > 1$ is due to \cite{Wiegner84}. 

The basic idea to prove
holomorphy of $E_j(c^{-1})$ w.r.t.\  $c^{ - 2}$ is the following: 
Since $(H (c) - m c ^{2} - z)^{ -1}$ is normal for
$z \in \C \setminus \R$,   \eqref{210} implies that the projection $P_j(c^{ -1})$ 
onto the eigenspace of the eigenvalue
$(E_j(c^{- 1} ) - z)^{-1}$ is holomorphic w.r.t.\  $c^{- 1}$ near $c^{- 1}=0$. 
To prove that $ E_j(c^{- 1} )$ is
actually holomorphic w.r.t.\  $c^{- 2}$ near $c^{- 2}=0$ one calculates
\begin{align}
& \tilde P_j(c^{-2}) := B(c)   P_j(c^{-1}) B(c)^{-1}=
 \begin{pmatrix}
1 & 0 \\
0 & c 
 \end{pmatrix} 
 \begin{pmatrix}
p_j & 0 \\
0 & 0 
 \end{pmatrix} 
 \begin{pmatrix}
1 & 0 \\
0 & 1/c 
 \end{pmatrix}  \nonumber \\
& + \left[ \text{terms holomorphic w.r.t. } c^{ -2} \right] 
= \left[  \text{terms holomorphic w.r.t. } c^{ -2} \right].
\end{align}

Here $\tilde P_j(c^{-2})$   and $p_j$ are the corresponding projections associated with \\
$ B(c) (H(c) - m c^{2} - z)^{-1} B(c)^{-1}$ and $(H_ + - z)^{-1} $ 
of dimension $m_j$, respectively. 
Thus, $|| \tilde P_j(c^{-2}) ||$ is bounded as $c^{ -2} \to 0$  and, hence, Butler's theorem (\cite{kato80}, p.\ 70) 
proves that
$\tilde P_j(c^{-2})$ and $(E_j (c^{- 1} ) - z)^{-1} $ are actually holomorphic
 w.r.t.\   $c^{ -2} $ near $c^{ -2} = 0$.

\section{Singular Weyl--Titchmarsh--Kodaira Theory}

The purpose of the remaining sections is to survey some highlights obtained by Fritz
for one-dimensional Schr\"odinger operators
\begin{equation}\label{eqschroe1d}
H =-\frac{d^2}{dx^2} + q, \qquad \text{in } \hr=L^2(a,b), \quad -\infty \le a < b \le \infty.
\end{equation}
In what follows we will denote by $\tau$ the
underlying differential operator and by $H$ an associated self-adjoint operator (as determined by separated boundary
conditions whenever necessary, i.e., when $\tau$ is in the limit circle case at one of the endpoints).
We will assume that the reader is familiar with the basic results from spectral theory
for these operators and refer, e.g., to the textbooks \cite{ls2,tschroe,wd,wd2}.

One key object in direct and inverse spectral theory is the Weyl--Titchmarsh $m$-function.
To define it suppose that the endpoint $a$ is regular, that is $a$ is finite and all solutions
(as well as their derivatives) extend continuously to this endpoint. For simplicity take
a Dirichlet boundary condition, $f(a)=0$, at $a$. Then we can choose a basis of solutions
of $\tau u = z u$ satisfying the initial conditions
\be\label{tpic}
\theta(z,a)=\phi'(z,a)=1, \qquad \theta'(z,a)=\phi(z,a)=0.
\ee
In particular, $\phi(z,x)$ satisfies the Dirichlet boundary condition at $a$. Moreover, for $z\in\C\setminus\R$ there will also be a
unique (up to scaling) solution $\psi(z,x)$ which will be square integrable near the other end point $b$
and which satisfies a possible boundary condition at $b$. If we normalize this function according to
$\psi(z,a)=1$ it can be written as the following linear combination of our basis of solutions
\be\label{defm}
\psi(z,x) = \theta(z,x) + m(z) \phi(z,x),
\ee
where the Weyl--Titchmarsh $m$-functions is given by
\be
m(z) = \frac{\psi'(z,a)}{\psi(z,a)}.
\ee
It is a key result that $m(z)$ is a Herglotz--Nevanlinna function (i.e., an analytic function mapping the upper half plane into
itself) and hence it can be represented as
\be
m(z) = \re(m(\I)) + \int_\R \Big(\frac{1}{\lam-z} - \frac{\lam}{1+\lam^2}\Big) d\rho(\lam),
\ee
where $\rho$ is a Borel measure with $\int_\R (1+\lam^2)^{-1} d\rho(\lam)<\infty$ given by the Stieltjes inversion formula
\be\label{defrho}
\frac{1}{2} \left( \rho\big((\lam_0,\lam_1)\big) + \rho\big([\lam_0,\lam_1]\big) \right)=
\lim_{\eps\downarrow 0} \frac{1}{\pi} \int_{\lam_0}^{\lam_1} \im\big(m(\lam+\I\eps)\big) d\lam.
\ee
Moreover, the map
\be\label{defft}
f(x) \mapsto \hat{f}(\lam) = \int_a^b \phi(\lam,x) f(x) dx,
\ee
initially defined on functions $f$ with compact support in $[a,b)$,
extends to a unitary map from $L^2(a,b)$ to $L^2(\R,d\rho)$ which maps our operator $H$
to multiplication with the independent variable $\lam$. In particular, the spectral measure $\rho$
contains all the spectral information of $H$! In his paper with Zinchenko \cite{gzi}, Fritz
not only gave a particularly simple proof of this fact but also pointed out a crucial extension
which is now known as singular Weyl--Titchmarsh--Kodaira theory.

In fact, the assumption that one endpoint is regular excludes many important examples which
have a simple spectrum and where such a transformation should still exist by the spectral
theorem. Two examples of central interest in quantum mechanics are the
radial part of the Coulomb Hamiltonian, where $q(x)=\frac{l(l+1)}{x^2} -\frac{\gam}{x}$ on $(0,\infty)$, or
the harmonic oscillator, $q(x)=x^2$ on $(-\infty,\infty)$. In both cases it is still possible to
define an entire (w.r.t.\ $z$) solution $\phi(z,x)$ which is square integrable near the left endpoint $a$ and
satisfies a possible boundary condition there. Moreover, a singular Weyl $m$-function can
be defined as before once a second linearly independent solution $\theta(z,x)$ is chosen.
In fact, this observation was already made by Kodaira \cite{ko}. However, as pointed out by Kac \cite{ka},
this is nontrivial if $\theta(z,x)$ is required to be entire as well. In fact, we can always choose a second
entire solution but the Wronskian $W(\theta(z),\phi(z))=\theta(z,x)\phi'(z,x)-\theta'(z,x)\phi(z,x)$ might have
zeros! In \cite{gzi} this problem was overcome by showing that to define the spectral measure via \eqref{defrho} it
suffices if $\theta(z,x)$ is analytic in a neighborhood of the real line. And it is easy
to check that if $\theta_c(z,x)$, $\phi_c(z,x)$ are two entire solutions defined by the initial conditions
 $\theta_c(z,c)=\phi_c'(z,c)=1$ and $\theta_c'(z,c)=\phi_c(z,c)=0$, then
 \be\label{deftheta}
 \theta(z,x)= \frac{\phi'(z,c)}{\phi(z,c)^2+\phi'(z,c)^2} \theta_c(z,x) + \frac{\phi(z,c)}{\phi(z,c)^2+\phi'(z,c)^2} \phi_c(z,x)
 \ee
will do the trick. In fact, by construction it is meromorphic in $\C$ with all poles away from the real axis and
one easily checks $W(\theta(z),\phi(z))=1$ (evaluate the Wronskian at $x=c$ and recall that it
is independent of $x$).

\begin{theorem}[\cite{gzi}]
Suppose $H$ has an entire solution $\phi(z,x)$ which is square integrable near the left endpoint $a$ and
satisfies a possible boundary condition there. Let $\theta(z,x)$ be defined as in \eqref{deftheta}.

Then there exists a Weyl $m$-function which is meromorphic in $\C\setminus\R$ (without any poles accumulating near the real line)
defined via \eqref{defm} and an associated spectral measure defined via \eqref{defrho}. Moreover,
the map \eqref{defft} initially defined on functions $f$ with compact support in $[a,b)$,
extends to a unitary map from $L^2(a,b)$ to $L^2(\R,d\rho)$ which maps our operator $H$
to multiplication with the independent variable $\lam$.
\end{theorem}

The elegant proof from \cite{gzi} alluded to before now simply involves the following two well-known facts: First,
the resolvent of $H$ which is given by
\be
(H-z)^{-1} f(x) = \int_a^b G(z,x,y) f(y) dy,
\ee
where
\be\label{defgf}
G(z,x,y) = \begin{cases} \phi(z,x) \psi(z,y), & y\ge x,\\
\phi(z,y) \psi(z,x), & y\le x.\end{cases}
\ee
And second, Stone's formula
\be
\spr{f}{F(H) f} = \frac{1}{\pi} \lim_{\eps\downarrow 0} \int_\R F(\lam) \im\spr{f}{(H-\lam-\I\eps)^{-1}f}) d\lam
\ee
for any bounded continuous function $F\in C_b(\R)$. Now insert the definitions and compute
the limit to obtain
\be
\spr{f}{F(H) f} = \int_\R F(\lam) |\hat{f}(\lam)|^2 d\rho(\lam),
\ee
which in turn implies the theorem.

This seminal result initiated the development of Weyl--Titchmarsh theory at a singular endpoint and
has triggered a number of results recently \cite{E12,EGNT12,EGNT12a,ET12a,ET12,ful08,fl,FLL12,KST12a,KST12,KT11,KT12,KL11}.

\section{Inverse Spectral Theory and Trace Formulas}
\label{secIST}

Next we turn to trace formulas, another area where Fritz has made profound contributions. To this end suppose that $(a,b)=\R$ and
that the potential $q(x)$ is periodic: $q(x+1)=q(x)$. Then it is well known by Floquet theory that
the spectrum consists of an infinite number of bands
\begin{equation}
\sig(H) = [E_0,E_1] \cup [E_2,E_3] \cup \cdots,
\end{equation}
where $E_0<E_1 \le E_2 < E_3 \le \cdots$. Moreover, if we restrict our operator to the interval $(x,x+1)$
and impose Dirichlet boundary conditions $f(x)=f(x+1)=0$ at the endpoints, then we obtain a
sequence of eigenvalues $\mu_1(x)<\mu_2(x)<\cdots$ depending on the base point $x$. Again
it follows from Floquet theory that the Dirichlet eigenvalues $\mu_j(x)$ lie in the closures of the spectral
gaps:  $E_{2j-1} \le \mu_j(x) \le E_{2j}$.

Then the following trace formula
\be\label{trfperiodic}
q(x) = E_0 + \sum_{j=1}^\infty \left[ E_{2j-1} + E_{2j} - 2\mu_j(x) \right]
\ee
holds under some suitable assumptions on the potential (e.g., $q\in C^1$). This trace formula was first
obtained in the case where the spectrum has only finitely many gaps (and thus the above sum is finite)
and later on generalized to periodic and even some classes of almost-periodic potentials. Furthermore,
trace formulas were also known in the context of scattering theory, but in this case the right-hand side contains both
a sum corresponding to the eigenvalues plus an integral corresponding to the continuous spectrum.
Hence, Fritz was asking for a generalization of \eqref{trfperiodic} to more general potentials.

To this end let us briefly sketch one way of proving \eqref{trfperiodic}. Denote by $\psi_\pm(z,x)$ the
Floquet solution of the underlying differential equation and recall that the diagonal Green function
is given by
\be
G(z,x,x)= \frac{\psi_+(z,x)\psi_-(z,x)}{W(\psi_+(z),\psi_-(z))},
\ee
where
\begin{equation}
W_x(f,g) = f(x)g'(x) - f'(x)g(x)
\end{equation}
denotes the Wronskian of two absolutely continuous functions. Recall that the Wronskian of two
solutions corresponding to the same spectral parameter $z$ is independent of $x$. Then,
using the fact that solutions of the underlying differential equation $\tau u = z u$ corresponding to
constant initial conditions are entire functions of order one half, one can obtain the following
product representation
\be\label{prodg}
G(z,x,x) = \frac{1}{2\sqrt{E_0-z}} \prod_{j=1}^\infty \frac{z-\mu_j(x)}{\sqrt{(E_{2j-1}-z)(E_{2j}-z)}}.
\ee
Comparing this product representation with the well-known asymptotics
\be\label{asymg}
G(z,x,x)= \frac{1}{2\sqrt{-z}} \left(1 + \frac{q(x)}{2 z} + o(z^{-1})\right)
\ee
establishes \eqref{trfperiodic}.

How can one possibly generalize this argument? Since \eqref{asymg} holds for general potentials,
a generalization of \eqref{prodg} needs to be found and, as observed by Fritz, the correct
starting point is the following exponential Herglotz representation
\be\label{gfxi}
G(z,x,x) = |G(\I,x,x)| \exp\left( \int_\R \Big(\frac{1}{\lam-z} - \frac{\lam}{1+\lam^2}\Big) \xi(\lam,x) d\lam\right).
\ee
Here the $\xi$-function is given by the Stieltjes inversion formula
\be\label{defxi}
\xi(\lam,x) = \lim_{\eps\downarrow 0} \frac{1}{\pi}\arg G(\lam+\I\eps,x,x),
\ee
where the limit exists for a.e.\ $\lam\in\R$. Note that since $G(.,x,x)$ is a Herglotz--Nevanlinna function, we have $0\le \xi(\lam,x) \le 1$.

Now a trace formula which works for arbitrary potentials follows by comparing asymptotics as
before. This is one of the key results obtained together with Barry Simon in \cite{gsxi}.

\begin{theorem}[\cite{gsxi}]
Suppose $q$ is continuous and bounded from below. Choose $E_0 \le \inf\sig(H)$. Then
\be\label{trfxi}
V(x) = E_0+ \lim_{\eps\downarrow 0} \int_{E_0}^\infty \E^{-\eps\lam} (1-2\xi(\lam,x)) d\lam.
\ee
\end{theorem}

In fact, even more was shown in \cite{gsxi}. Namely, let $H^D_x = H^D_{x,-}\oplus H^D_{x,+}$ be the operator $H$ restricted according to
$L^2(a,b)= L^2(a,c) \oplus L^2(c,b)$ by imposing Dirichlet boundary conditions at the point $x\in(a,b)$. Then
$\xi(\lam,x)$ was identified as the Krein spectral shift function of the pair $(H,H^D_x)$, which in turn opened the door
for deriving numerous novel trace formulas for Schr\"odinger operators: \cite{g,gh1,gh2,gh3,gh4,ghs,ghsz1,ghsz2,ghsz3}.

To demonstrate the usefulness and generality of this theorem let us extract two special cases. First let us see how to
obtain \eqref{trfperiodic} by evaluating \eqref{defxi} using \eqref{prodg} (note that the limit is not needed since $G(\lam,x,x)$
has a continuous extension to the real line away from the band edges). We begin with $\lam<E_0$ in which case $G(\lam,x,x)$
is positive implying $\xi(\lam,x)=0$ for $\lam\in(-\infty,E_0)$. At $E_0$ the Green function $G(\lam,x,x)$ will become purely
imaginary implying $\xi(\lam,x)=\frac{1}{2}$ for $\lam\in(E_0,E_1)$. After $E_1$ it will become negative and change sign at its
zero $\mu_1(x)$, implying $\xi(\lam,x)=1$ for $\lam\in(E_1,\mu_1(x))$ and $\xi(\lam,x)=0$ for $\lam\in(\mu_1(x),E_2)$.
Clearly this pattern keeps repeating and we obtain \eqref{trfperiodic} (assuming $q\in C^1$ the spectral gaps will close
sufficiently fast, such that the limit in \eqref{trfxi} can be taken inside the integral).

A novel application is to confining potentials, like the harmonic oscillator, satisfying $q(x)\to +\infty$ as $|x| \to\infty$.
Then $H$ has purely discrete spectrum $E_0<E_1<\cdots$ and the same is true for $H^D_x$ whose eigenvalues
$\mu_1(x)<\mu_2(x)<\cdots$ are known to satisfy $E_{j-1} \le \mu_j(x) \le E_j$. Then note that $G(z,x,x)$ has first
order poles at the eigenvalues $E_j$ and first order zeros at the Dirichlet eigenvalues $\mu_j(x)$ (since either
$\psi_-(\lam,x)$ or $\psi_+(\lam,x)$ must vanish at $\lam=\mu_j(x)$ --- if both should vanish they must be linearly
dependent and thus such a point must also be an eigenvalue). Hence $G(z,x,x)$ will be real-valued on the real
line and change sign at every pole $E_j$ and every zero $\mu_j(x)$ implying
\be
V(x) =  E_0+ \lim_{\eps\downarrow 0} \frac{1}{\eps} \sum_{j\in\N} 2\E^{-\eps \mu_j(x)} -\E^{-\eps E_{j-1}} - \E^{-\eps E_j},
\ee
which is just an abelianized version of \eqref{trfperiodic}.

This brings us to another item discussed in \cite{gsxi}, namely the consequences of the above ideas for inverse spectral theory.
While the inverse spectral theory is well understood for periodic operators, this is not the case for confining potentials.
For example, one open question is to describe the isospectral class of the harmonic oscillator. To see how the
$\xi$-function can help understanding this problem, we begin with the observation that $\xi(\lam,x)$ determines $G(z,x,x)$
by virtue of \eqref{gfxi} (the unknown constant can be determined from the asymptotics \eqref{asymg}). 

Now let $\psi_\pm(z,x)$ be the solutions of $\tau u = z u$ which are square integrable near $\pm\infty$ and observe that
\eqref{defgf} reads
\be
G(z,x,x)= \frac{-1}{m_+(z,x)+m_-(z,x)},
\ee
where
\be
m_\pm(z,x) =\pm \frac{\psi_\pm'(z,x)}{\psi_\pm(z,x)}
\ee
are the Weyl $m$-functions of $H^D_{x,\pm}$. Hence if you fix $x$, say $x=0$, then $\xi(\lam,0)$ determines the
sum of the spectral measures $\rho_{0,\pm}$ corresponding to $H^D_{0,\pm}$ and it remains to split this
information into its two pieces (recall that $\rho_{0,\pm}$ determine $H^D_{0,\pm}$ and thus $H$).
Clearly this is impossible in general without additional data. However, it turns out that things get particularly simple
if the following reflectionless property is assumed:
\be
m_-(\lam+\I0,x) = - \ol{m_-(\lam+\I0,x)}, \qquad \text{a.e. } \lam\in\sig_{ess}(H),
\ee
for one (it then follows for all) $x$. This covers for example periodic operators, reflectionless potentials from
scattering theory (which is where the name comes from), or operators with purely discrete spectrum.
In fact, this condition tells us that the absolutely continuous part of the measure needs to be split equally.
Moreover, assuming that there is no singularly continuous part it remains to assign a sign $\sig_j\in\{\pm\}$
to every Dirichlet eigenvalue $\mu_j=\mu_j(0)$ which encodes if it is an eigenvalue of $H^D_{0,-}$ or
$H^D_{0,+}$ (if it is an eigenvalue of both, which can only happen if it is also an eigenvalue of $H$
as pointed out before, a number $\sig_j\in (-1,1)$ is needed to encode how the mass should be split).

By construction we obtain that the eigenvalues $\{ E_j \}_{j\in\N_0}$ together with the Dirichlet data
$\{(\mu_j,\sig_j)\}_{j\in\N}$ uniquely determine $H$. The remaining problem to determine the isospectral
class of a given operator $H_0$ is to describe the set of admissible Dirichlet data. One key contribution
towards answering this question is the Dirichlet deformation method to be discusses in the next
section, which shows that all restrictions (apart form the obvious ones already listed above) must
be of an asymptotic nature since any finite part of the Dirichlet data can be changed at will.

For further information on these circle of ideas see \cite{cg,cgr,g,gh1,gh2,gh3,gh4,ghs,ghsz1,ghsz2,ghsz3,gmt,grt,gsp,ps,ry1,ry2,ry3,sh,t6,zi}.
In particular, we mention also the review by Fritz \cite{gefest}.

\section{Commutation Methods} \label{commutation}

Another topic where Fritz made seminal contributions is commutation methods. These are
methods of inserting (and removing) eigenvalues in spectral gaps of a given
one-dimensional Schr\"odinger operator $H$. They play a prominent role in diverse fields such as the inverse scattering approach
(see, e.g., \cite{dt}, \cite{mc} and the references therein), supersymmetric quantum mechanics
(see, e.g., \cite{gss} and the references therein), level comparison theorems
(cf$.$ \cite{ba} and the literature cited therein), and in
connection with B\"acklund transformations for the KdV hierarchy (see, e.g., \cite{gs}, \cite{gss} and the references therein).

Historically, these methods of inserting eigenvalues go back to Jacobi
\cite{ja} and Darboux \cite{da} with decisive later contributions by Crum
\cite{cr}, Schmincke \cite{sc}, and, especially, Deift \cite{de}. Two
particular methods, shortly to be discussed in an informal manner in
\eqref{tau}--\eqref{qdc} below, turned out to be of special importance: The
single commutation method, also called the Crum--Darboux method \cite{cr},
\cite{da} (actually going back at least to Jacobi \cite{ja}) and the double
commutation method, to be found, e.g$.$, in the seminal work of Gel'fand and
Levitan \cite{gl}.

The single commutation method relies on existence of a positive solution $\psi$
of $H \psi = \lam \psi$ which confines its applicability to the insertion of
eigenvalues below the spectrum of $H$ (assuming $H$ to be bounded from below).
Introducing
\begin{equation}\label{tau}
A = \frac{d}{dx} + \phi,
\quad A^* =  -\frac{d}{dx} + \phi, \qquad\text{where } \phi = \frac{d}{dx} \log\psi,
\end{equation}
a straightforward calculation reveals
\begin{equation}
H = A^* A +\lam= -\frac{d^2}{dx^2} +q , \quad \hat{H} = A A^* +\lam = - \frac{d^2}{dx^2}+ \hat{q},
\end{equation}
with
\begin{equation} \label{qsc}
\hat{q} = q - 2 \frac{d^2}{dx^2} \log\psi.
\end{equation}
Thus (taking proper domain considerations into account) we can define two
operators $H,\hat{H}$ on $\hr=L^2(a,b)$ which turn out unitarily equivalent when restricted to the orthogonal complement of
the eigenspaces corresponding to $\lam$. In fact, using the polar decomposition $A = U |A|$, where
$|A|= (A^*A)^{1/2}$ and $U : \ker(A)^\perp \to \ker(A^*)^\perp$ is unitary, one infers
$\hat{H} P = (A A^* +\lam ) P = U |A| |A| U^* + \lam P = U HU^*$ with $P=UU^*$, the projection onto $\ker(A^*)^\perp$.
Moreover, $H-\lam, \hat{H}-\lam\ge 0$ which is equivalent to the existence of
the positive solution $\psi$ \cite{gz}. Formulas \eqref{tau}--\eqref{qsc}
constitute the single commutation method.

The double commutation method on the other hand, allows one to insert
eigenvalues into {\em any} spectral gap of $H$. To this end we assume that
$\psi$ is square integrable near $a$ and  consider two more expressions
$A_\gam$, $A_\gam^*$ as above with $\psi_\gam = \psi/(1+\gam \int_a^x
\psi(t)^2 dt)$. This implies
\begin{equation} \label{dc}
H= A A^*+\lam  = A_\gam^* A_\gam +\lam, \quad H_\gam = A_\gam A_\gam^*
+ \lam = - \frac{d^2}{dx^2} + q_\gam,
\end{equation}
where
\begin{equation} \label{qdc}
q_\gam = q - 2\frac{d^2}{dx^2} \log
\Big(1+\gam \int_a^x \psi(t)^2 dt \Big).
\end{equation}
The considerations for the single commutation method show that $H$ and $H_\gam$ are
unitarily equivalent when restricted to the orthogonal complement of the eigenspaces corresponding to $\lam$
as long as $\psi$ is positive. However, observe that $q_\gam$ is well defined even if $\psi$ has zeros, and
it is natural to conjecture that the last claim continuous to hold even in the case where all intermediate operators are
ill-defined. This turned out much harder to prove and was achieved by Fritz
\cite{fg} on the basis of Weyl--Titch\-marsh $m$-function techniques. 

In what follows we will denote by $\tau$ the
underlying differential operator and by $H$ an associated self-adjoint operator (as determined by separated boundary
conditions whenever necessary, i.e., when $\tau$ is in the limit circle case at one of the endpoints).

To this end, suppose $a$ is regular (i.e., $a$ is finite
and $q\in L^1(a,c)$ such that all solutions extend continuously to the endpoint $a$) and introduce the 
Weyl--Titch\-marsh $m$-function as
\begin{equation}
m(z)= \frac{\psi'(z,a)}{\psi(z,a)},
\end{equation}
where $\psi(z,x)$ is the unique (up to scaling) solution of $\tau\psi=z \psi$ which is square integrable and
satisfies a possible boundary condition at $b$. Then it is well known that $m(z)$ carries all the information on $H$
and hence it suffices to find an explicit expression for the Weyl--Titch\-marsh $m$-function of $H_\gam$ in terms
of the one of $H$. This strategy was carried out in \cite{fg} and the following theorem was obtained as one of the
main results.

\begin{theorem}[\cite{fg,gtdc}]
Let $\lam \in \R$ and $\psi(\lam,.)$ be a solution satisfying the following conditions:
\begin{enumerate}
\item $\psi, \psi' \in AC_{loc}(a,b)$ and $\psi$
is a real-valued solution of $\tau \psi = \lam \psi$.
\item $\psi$ is square integrable near $a$ and fulfills the
boundary condition (of $H$) at $a$ and $b$ if any (i.e$.$, if $\tau$ is limit
circle ($l.c.$) at $a$ respectively $b$).
\end{enumerate}
Define
\begin{equation}
\psi_\gam(\lam,x) = \frac{\psi(\lam,x)}{1 + \gam \int_a^x \psi(\lam,t)^2 dt}.
\end{equation}
Let $P(\lam)$, $P_\gam(\lam)$ be the projections on the subspaces
spanned by $\psi(\lam)$, $\psi_\gam(\lam)$, respectively (if one of these
functions is not square integrable, set the corresponding projection equal to zero).
Then the operator $H_\gam$ defined by
\begin{equation}
H_\gam f = \tau_\gam f, \quad \db(H_\gam) = \bay[t]{l} \{ f \in 
\hr |\, f,pf' \in AC_{loc}((a,b)) ; \tau_\gam f \in \hr  ;\\
W_a(\psi_\gam(\lam),f)=W_b(\psi_\gam(\lam),f)=0 \}, \eay
\end{equation}
with $q_\gam$ given by \eqref{qdc} is self-adjoint. Moreover, $H_\gam$ has the eigenfunction
$\psi_\gam(\lam)$ associated with the eigenvalue $\lam$. If $\psi(\lam) \not\in \hr$ (and hence
$\tau$ is limit point ($l.p.$) at $b$) we have
\begin{equation}
 H_\gam (1-P_\gam(\lam)) = U_\gam H U_\gam^{-1} (1-P_\gam(\lam))
\end{equation}
for some unitary operator $U_\gam$ and thus
\begin{equation}
\bay{rcl@{\qquad}rcl}
\sigma(H_\gam) &=& \sigma(H) \cup \{ \lam\}, & \sigma_{ac}(H_\gam)
&=& \sigma_{ ac}(H), \\  \sigma_{p}(H_\gam) &=& \sigma_{p}(H) \cup \{
\lam\}, & \sigma_{sc}(H_\gam) &=& \sigma_{sc}(H).
\eay 
\end{equation}
(Here $\sigma_{ac}(.),\sigma_{sc}(.)$ denotes the absolutely and singularly
continuous spectrum, respectively.) If $\psi(\lam) \in \hr $ there
is a unitary operator $\ti{U}_\gam = U_\gam \oplus
\sqrt{1+\gam \| \psi(\lam) \|^2} U_\gam$ on $(1-P(\lam)) \hr 
\oplus P(\lam) \hr $ such that
\begin{equation}
H_\gam = \ti{U}_\gam H \ti{U}_\gam^{-1}.
\end{equation}
\end{theorem}

As already mentioned, the main part of this result is due to Fritz \cite{fg}. In Gesztesy and Teschl \cite{gtdc}
the transformation operator was identified to be
\begin{equation}
(U_\gam f)(x) = f(x) - \gam \psi_\gam(\lam,x) \int_a^x \psi(\lam,t) f(t) dt,
\end{equation}
and the above theorem reduces to checking that $U_\gam$ has the claimed properties. Moreover, the
whole method was extended to Sturm--Liouville operators and some technical assumptions were relaxed.
For the connections with singular Weyl--Titchmarsh--Kodaira theory see \cite{KST12a}.

A further decisive contribution was made by Gesztesy, Simon, and Teschl \cite{gstddm} motivated by a commutation method
first introduced by Finkel, Isaacson, and Trubowitz \cite{fit} in connection with an explicit realization of the isospectral
torus of periodic potentials. This method was again used by Buys and Finkel \cite{bufi} (see also Iwasaki \cite{iwa})
in the context of periodic finite-gap potentials and by P\"oschel and Trubowitz \cite{potr} and
Ralston and Trubowitz \cite{ratr} for various boundary value problems on compact intervals.
As in the previous case, this method formally consists of two single commutations, but this
time at different values of the spectral parameter $\mu$ and $\nu$. The resulting operator will have the Wronskian $W(\psi(\mu),\psi(\nu))$ in its denominator
and it turns out that it will be well-defined as long as both $\mu$ and $\nu$ lie in the same spectral gap of $H$. 
This is related to the fact, that $(H-\mu)(H-\nu)$ will still be nonnegative
under this assumption and hence can still be factorized as $(H-\mu)(H-\nu)=B^*B$, as explained later by Schmincke \cite{sc2}.
Moreover, the fact that the zeros of the Wronskian of two solutions is related to the spectrum lead to the development of
renormalized oscillation theory \cite{gstz} to be discussed in Section~\ref{secOT} below.

To explain the relevance of this method recall that an important role in inverse spectral theory
is played by the so called Dirichlet data as introduced at the end of the previous section. The key result of
\cite{gstddm} states that the above commutation method can be used to move a given Dirichlet eigenvalue
to any other admissible position within its gap as well as change its {\em sign}.  Hence it is also known as
the Dirichlet deformation method. We refrain from further
details here and refer to Fritz's own review \cite{gefest} instead. For generalizations to other operators see
\cite{sa0,sa1,t2,t3,tjac}.

For further generalizations of commutation methods see \cite{bbw,bt,ahm}. Finally, these methods are also
relevant in connection with inverse scattering theory
(see, e.g., \cite{ak,mc} and the references cited therein)
and yield a direct construction of $N$-soliton solutions relative to arbitrary
background solutions of the (generalized) KdV hierarchy along the
methods of \cite{gs} (see also \cite{gu}).
Generalizations to other integrable equations can be found in \cite{cghl,ou,bght,gr,sa2,t5}.

\section{Oscillation Theory}
\label{secOT}

Another area where Fritz had an important impact is oscillation theory. As indicated at the end of the previous section, a crucial
ingredient in the development of the Dirichlet deformation method was the fact that the Wronskian of two solutions
$W(\psi(\mu),\psi(\nu))$ is nonzero as long as there is no part of the spectrum inside the interval $(\mu,\nu)$. This observation
naturally lead to the question how the number of points in the spectrum between $\mu$ and $\nu$ are related to the
number of zeros of the Wronskian. This question was answered by Gesztesy, Simon, and Teschl in \cite{gstz}.

Denote by $W_0(f,g)$ be the number of zeros of the Wronskian in the open interval $(a,b)$ not counting multiplicities of zeros. Given
$E_1<E_2$, we let $N_0(E_1, E_2)=\dim\ran\, P_{(E_1, E_2)}(H)$ be the dimension of the spectral
projection $P_{(E_1, E_2)}(H)$ of $H$. Then the main result from \cite{gstz} is the following:
 
\begin{theorem}[\cite{gstz}]
Let $\psi_-(\lam,x)$,  $\psi_+(\lam,x)$ be solutions of $\tau \psi = \lam \psi$ which are square integrable near $a$, $b$ and
satisfy a possible boundary condition at $a$, $b$, respectively. (Such solutions will in general only exist inside a spectral gap).
Suppose $E_1 < E_2$. Then
\be
W_0(\psi_-(E_1),\psi_+(E_2))=N_0(E_1, E_2).
\ee
\end{theorem}

This result constituted an important generalization suitable for counting the number of eigenvalues
inside essential gaps. In fact, there has been considerable efforts to generalize classical oscillation theory
to singular operators by Hartman and others around 1950. However, while these efforts were successful below the
essential spectrum, they were only partly successful for counting eigenvalues in essential spectral gaps. For example
Hartmann \cite{har} could show that for a regular endpoint $a$ and a limit point endpoint $b$, counting the difference of zeros of the solutions $\psi_-(E_2,x)$ and
$\psi_-(E_1,x)$ on $(a,c)$ will have $N_0(E_1, E_2)$ as its $\liminf$ when $c\to b$. As a simple example, the above
theorem allows to extend this result to the case when $a$ is non-oscillatory (thus covering important examples like radial
Schr\"odinger equations).

Finally let us outline a proof for the above result. The first step is to prove this result in the case where both endpoints are
regular. This can be done using the usual Pr\"ufer techniques and we refer to the textbook \cite{tode} for a simple proof.
Next one can approximate $H$ on $(a,b)$ by regular operators on $(c,d)\subseteq (a,b)$. Again this technique is
standard (for a nice review see \cite{wdapprox}) and the restricted operators will converge to $H$ in the strong resolvent sense. Since the spectrum
cannot expand for such limits, one obtains $W_0 \ge N_0$. To obtain the reverse inequality one uses a variational argument.
To this end observe that at every zero of the Wronskian the two solutions $\psi_-(E_1)$ and $\psi_+(E_2)$ can be glued together
(by scaling one to match the value of the other) to give a trial function in the domain of $H$ (since the derivatives will then automatically
match at this point). Hence one obtains a $W_0$ dimensional space of trial functions $\psi$ satisfying
\be
\big\| \big(H - \frac{E_2-E_1}{2}\big) \psi\big\| < \frac{E_2-E_1}{2} \|\psi\|
\ee
and the required reverse inequality $W_0 \le N_0$ follows from the spectral theorem. For further recent reviews of these methods
we refer to the ones by Fritz \cite{gefest} and Simon \cite{siosc}.

Again this result has triggered several extensions. In particular, it was shown by Kr\"uger and Teschl \cite{kt} that one can take
solutions of different operators if the right-hand side is interpreted as spectral shift between these two operators.
We refer to \cite{kt,kt2} and the references therein. The question whether eigenvalues accumulate at the boundary
of an essential spectral gap based on these methods is considered in \cite{kt3,kms}. In this respect we should also
mention the beautiful result by Fritz and \"Unal \cite{gun} which gives the most general version of Kneser's theorem.
Extensions to other operators can be found in \cite{at,el,st,t1,t4}.

\section{Non-self-adjoint operators}

While all results so far were concerned with self-adjoint operators, Fritz always was quite active in the area of non-self-adjoint operators
as well. As a prototypical example we mention his beautiful work with Tkachenko \cite{gta1,gta2} solving the long-standing open question when
a non-self-adjoint Hill operator is a spectral operator of scalar type in the sense of Dunford. This problem had been open for about 40 years. 

We have already encountered Hill's equation (i.e, \eqref{eqschroe1d} with periodic potential $q(x+1)=q(x)$) in Section~\ref{secIST}.
But now we will allow $q$ to be complex valued (in addition it is assumed to be locally square integrable).
Recall that one of the key objects from Floquet theory is the Floquet discriminant
\be
\Delta(z)= \frac{\theta(z,1)+\phi'(z,1)}{2},
\ee
where $\theta(z,x)$, $\phi(z,x)$ is a fundamental system of solutions satisfying the initial conditions \eqref{tpic} at $a=0$.
It was shown by Serov \cite{se} that the spectrum of $H$ is given by $\sig(H) =\{ z \in \C | \Delta(z)\in [-1,1] \}$. In particular,
the spectrum consists of a number of analytic arcs which, however, might intersect in inner points \cite{pa}.

The following version of this criterion involves the spectrum $\sigma(H)$ of $H$, the
Dirichlet spectrum $\{\mu_k\}_{k\in\N}$, the periodic
spectrum $\{E_k(0)\}_{k\in {\N_0}}$, and the antiperiodic spectrum
$\{E_k(\pi)\}_{k\in {\N_0}}$, and is connected with the algebraic
and geometric multiplicities of the eigenvalues in the sets
$\sigma(H(t))$, where $H(t)$, $t\in [0, 2\pi]$, denotes the densely defined closed realization of \eqref{eqschroe1d} in $L^2(0,1)$
in terms of the $t$-dependent boundary conditions $f(1)=\E^{\I t}f(0)$, $f'(1)=\E^{\I t}f'(0)$,  $t\in [0, 2\pi]$. 
It is known that the union of the periodic and antiperiodic spectra is formed by the numbers
\begin{align}
& \lambda_k^\pm=\left(k+\frac{\int_0^1
q(x) dx}{2k}+\frac{s_k^\pm}{k}\right)^2,
\; k\in\N, \quad
\sum_{k\in{\mathbb\N}}|s_k^\pm|^2<\infty,  \\
& \{E_k(0)\}_{k\in\N_0}
= \{\lambda_0^+, \lambda_{2k}^+, \, \lambda_{2k}^-\}_{k\in\N},    \quad
 \{E_k(\pi)\}_{k\in\N_0}
= \{\lambda_{2k+1}^+, \, \lambda_{2k+1}^-\}_{k\in\N_0}, 
\end{align}

\begin{theorem}\label{tGT1}
A Hill operator $H$ is a spectral operator of scalar type if and only if
the following conditions $(i)$ and $(ii)$ are satisfied: \\
$(i)$ For all $t\in[0,2\pi]$ and all $E_k(t)\in\sigma(H(t))$,
each root function of the operator $H(t)$ associated with $E_k(t)$ is an
eigenfunction of $H(t)$.  \\
$(ii)$ Let
\begin{equation}\label{l610}
\mathcal{Q}=\{k\in\N\,|\, d_k={\rm dist}(\delta_k,\sigma(H))>0\}.
\end{equation}
Then
\begin{equation} \label{l60}
\displaystyle
\sup_{k\in\mathcal{Q}}\frac{|\lambda_k^+-\lambda_k^-|}
{{\rm dist}(\delta_k,\sigma(H))}<\infty,\quad
\sup_{k\in\mathcal{Q}}\frac{|\mu_k-\lambda_k^-|}
{{\rm dist}(\delta_k,\sigma(H))}<\infty,\quad
\sup_{k\in\mathcal{Q}}\frac{|\mu_k-\lambda_k^+|}
{{\rm dist}(\delta_k,\sigma(H))}<\infty.
\end{equation}
\end{theorem}

Here a root function of $H(t)$ associated with the eigenvalue $E_k(t)$
denotes any element $\psi$ satisfying $(H(t)-E_k(t))^m\psi=0$ for some $m\in\N$ 
(i.e., any element in the algebraic eigenspace of $H(t)$ corresponding to $E_k(t)$). Of course, the conditions $(i)$ and $(ii)$ are always satisfied in the self-adjoint context where $q$ is real-valued. 

It can be shown that the conditions $(i)$ and $(ii)$ imposed in Theorem \ref{tGT1} yield the remarkable fact that the property of a Hill operator being a spectral operator is independent of smoothness (or even analyticity) properties of the potential $q$. In addition, also a functional model for periodic Schr\"odinger operators  that are spectral operators of scalar type was established and the corresponding eigenfunction expansion was developed in \cite{gta1}, \cite{gta2}. 

In their recent paper \cite{gta3}, under the assumption that $q \in L^2(0,1)$, necessary and sufficient conditions in terms of spectral data for (non-self-adjoint) Schr\"odinger operators in $L^2(0,1)$ with periodic and antiperiodic boundary conditions to possess a Riesz basis of root vectors were derived. Without entering details, we mention that this problem generated an enormous amount of interest and remained open for a long time. Their key result reads:

\begin{theorem} \label{tGT2}
Assume $q \in L^2(0,1)$, then the following results hold: \\ 
$(i)$ The system of root vectors of $H(0)$ contains a 
Riesz basis in $L^2(0,1)$ if and only if 
\begin{equation}
\sup_{\substack{k \in \N, \\ \lambda_{2k}^+ \neq \lambda_{2k}^-}} 
\frac{|\mu_{2k} - \lambda_{2k}^{\pm}|}{|\lambda_{2k}^+ - \lambda_{2k}^-|} < \infty.   
\label{1.19}
\end{equation}
$(ii)$ The system of root vectors of $H(\pi)$ contains a 
Riesz basis in $L^2(0,1)$ if and only if 
\begin{equation}
\sup_{\substack{k \in \N, \\ \lambda_{2k+1}^+ \neq \lambda_{2k+1}^-}} 
\frac{|\mu_{2k+1} - \lambda_{2k+1}^{\pm}|}{|\lambda_{2k+1}^+ - \lambda_{2k+1}^-|} < \infty.  \label{1.20}
\end{equation}
\end{theorem} 

Here $\sup_{k \in \N, \, \lambda_{j}^+ \neq \lambda_{j}^-}$ signifies that all subscripts $j\in\N$ in \eqref{1.19} and \eqref{1.20} for which $\lambda_{j}^+$ and $\lambda_{j}^-$ coincide are simply excluded from the supremum considered. 

One observes that only the simple periodic (resp., antiperiodic) eigenvalues enter in the necessary and sufficient conditions 
\eqref{1.19} (resp., \eqref{1.20}) for the existence of a Riesz basis of root vectors of $H(0)$ (resp., $H(\pi)$). The 
multiple periodic (resp., antiperiodic) eigenvalues play no role in deciding whether or not the system of root vectors of $H(0)$ (resp., $H(\pi)$) constitutes a Riesz basis in $L^2(0,1)$.
In addition, only every other Dirichlet eigenvalue (i.e., half the Dirichlet spectrum) enters the criterion \eqref{1.19} (resp., \eqref{1.20}). 

For additional detailed results in this direction we also refer to \cite{DM12} and the extensive literature cited therein.

\bigskip
\noindent
{\bf Acknowledgments.}
We thank Barry Simon for valuable comments improving the presentation of the material.

\end{document}